\documentclass[letterpaper, 10 pt, journal, twoside]{IEEEtran}
\usepackage[utf8]{inputenc}
\usepackage{amsmath}
\usepackage{amsfonts}
\usepackage{amssymb}
\let\labelindent\relax \usepackage{enumitem}  
\usepackage{color}
\usepackage{graphicx}
\usepackage[caption=false, font=footnotesize]{subfig}
\newtheorem{theorem}{Theorem}

\newtheorem{definition}{Definition}

\DeclareMathOperator{\tr}{trace}

\DeclareMathOperator{\sing}{sing}
\DeclareMathOperator{\mult}{mult}

\title{Willems' Fundamental Lemma for State-space Systems 
and its Extension to Multiple Datasets}

\author{Henk J. van Waarde, Claudio De Persis, M. Kanat Camlibel, and Pietro Tesi 
\thanks{The authors are with the Jan C. Willems Center for Systems and Control, 
University of Groningen, The Netherlands. 
Pietro Tesi is also with the Department of Information Engineering, University of Florence, Italy.
Email: {\tt\small \{h.j.van.waarde, c.de.persis, m.k.camlibel\}@rug.nl; pietro.tesi@unifi.it}.
}%
}

\begin{document}

\maketitle
\thispagestyle{empty}

\begin{abstract}
Willems \emph{et al.}'s fundamental lemma asserts that all trajectories of a linear system can be obtained from a single given one, assuming that a persistency of excitation and a controllability condition hold. This result has profound implications for system identification and data-driven control, and has seen a revival over the last few years. The purpose of this paper is to extend Willems' lemma to the situation where multiple (possibly short) system trajectories are given instead of a single long one. To this end, we introduce a notion of collective persistency of excitation. We will show that all trajectories of a linear system can be obtained from a given finite number of trajectories, as long as these are collectively persistently exciting. We will demonstrate that this result enables the identification of linear systems from data sets with missing samples. Additionally, we show that the result is of practical significance in data-driven control of unstable systems.  
\end{abstract}


\section{Introduction}

\IEEEPARstart{I}{n} the seminal work by Willems and coauthors \cite{Willems2005}, it was shown that a single, sufficiently exciting trajectory of a linear system can be used to parameterize \emph{all} trajectories that the system can produce. This result has later been named the \emph{fundamental lemma} \cite{Markovsky2005,Markovsky2008}, 
and plays an important role in the learning and control of dynamical systems on the basis of measured data. 

An immediate consequence of the fundamental lemma is that a 
persistently exciting trajectory 
captures the entire behavior of the 
data-generating system, thus allowing successful identification of a system model using subspace methods \cite{Verhaegen1992}. The lemma also enables data-driven simulation \cite{Markovsky2008}, which involves the computation of the system's response to a given reference input. In addition, Willems' lemma is instrumental in the design of controllers from data. The result has been applied to tackle several control problems, ranging from output matching \cite{Markovsky2008} to control by interconnection \cite{Maupong2017}, 
predictive control \cite{Coulson2019,Huang2019,Berberich2019b}, optimal and robust control \cite{DePersis2019}, linear quadratic regulation 
\cite{Markovsky2008,DePersis2019,Rotulo2019} 
as well as set-invariance control \cite{Bisoffi2019}.

All of the above examples show the value of the fundamental lemma in modeling, simulation and control using a \emph{single} measured system trajectory. Nonetheless, there are many scenarios in which \emph{multiple} system trajectories are measured instead of a single one. For example, performing multiple short experiments becomes desirable when the data-generating system has unstable dynamics. Also, as pointed out in \cite{Holcomb2017}, a single system trajectory collected during normal operations may be too poorly excited to reveal the system dynamics. In contrast, multiple archival data may \emph{collectively} provide a well-excited experiment. Another situation is when a single trajectory is measured but some of the samples are corrupted or missing. In this case, we have access to multiple system trajectories consisting of the remaining, uncorrupted, data samples. System identification from multiple experiments \cite{Markovsky2013b,Markovsky2015} and from data with missing samples \cite{Markovsky2013,Markovsky2016,Markovsky2017} has been studied. However, a proof of Willems' lemma for multiple trajectories is still missing. Therefore, in this paper we aim at extending Willems' fundamental lemma to the case where multiple trajectories, possibly of different lengths, are given instead of a single one.

Originally, the fundamental lemma was formulated and proven in a behavioral context. The starting point in this paper, however, is a reformulation of the lemma in terms of state-space systems. Such a version of Willems' fundamental lemma has appeared before in \cite[Lem. 2]{DePersis2019} and \cite[Thm. 3]{Berberich2019} but no proof of the statement was given in this context.
Our first contribution is to provide a complete and self-contained 
proof of the lemma for state-space systems.
Strictly speaking, such an alternative proof is not necessary since the original proof of \cite{Willems2005} applies to state-space systems as a special case. Nonetheless, we believe that our proof can be of interest to researchers who want to apply Willems' lemma to state-space systems. In fact, the proof is \emph{elementary} in the sense that it only makes use of basic concepts such as the Cayley-Hamilton theorem and Kalman controllability test. The proof is also direct, and in contrast to \cite{Willems2005} does not rely on a contradiction argument.

Our second contribution involves the extension of the fundamental lemma to the case of multiple trajectories. To this end, we first introduce a notion of \emph{collective} persistency of excitation. Then, analogous to Willems' lemma, we show that a finite number of given trajectories can be used to parameterize all trajectories of the system, assuming that collective persistency of excitation holds. We will illustrate this result by two examples. First, we will show that the extended fundamental lemma enables the identification of linear systems from data sets with missing samples. Next, we will show how the result can be used to compute controllers of unstable systems from multiple short system trajectories, even when this is problematic from a single long trajectory. 

The paper is organized as follows: in Section \ref{sectionWillemslemma} we formulate and prove Willems' fundamental lemma. Section \ref{sectionmultipledata} extends the lemma to multiple trajectories. In Section \ref{sectionexamples} we provide applications of this result. Finally, Section \ref{sectionconclusions} contains our conclusions. 

\subsection{Notation}

\noindent 
The \emph{left kernel} of a real matrix $M$ is the space of all real row vectors $v$ such that $vM = 0$. The zero vector of dimension $n$ is denoted by $0_n$. Consider a signal $f: \mathbb{Z} \to \mathbb{R}^\bullet$ and let $i,j \in \mathbb{Z}$ be integers such that $i \leq j$. We denote by $f_{[i,j]}$ the restriction of $f$ to the interval $[i,j]$, that is, 
\begin{equation*}
    f_{[i,j]} := \begin{bmatrix} f(i)^\top & f(i+1)^\top & \cdots & f(j)^\top
    \end{bmatrix}^\top.
\end{equation*}
With slight abuse of notation, we will also use the notation $f_{[i,j]}$ to refer to the sequence $f(i),f(i+1),\dots,f(j)$. 
Let $k$ be a positive integer such that $k \leq j-i+1$ and define the \emph{Hankel matrix} of depth $k$, associated with $f_{[i,j]}$, as
\begin{equation*}
    \mathcal{H}_{k}(f_{[i,j]}) := 
    \begingroup 
    \setlength\arraycolsep{2pt}
    \begin{bmatrix} f(i) & f(i+1) & \cdots & f(j-k+1) \\
    f(i+1) & f(i+2) & \cdots & f(j-k+2) \\
    \vdots & \vdots & & \vdots \\
    f(i+k-1) & f(i+k) & \cdots & f(j) 
    \end{bmatrix}.
    \endgroup
\end{equation*}
Note that the subscript $k$ refers to the number of block rows of the Hankel matrix. 
\begin{definition}
The sequence $f_{[i,j]}$ is said to be \emph{persistently exciting of order $k$} 
if $\mathcal{H}_{k}(f_{[i,j]})$ has full row rank. 
\end{definition}

\section{Willems \emph{et al.}'s fundamental lemma in the context of state-space systems}
\label{sectionWillemslemma}
In this section we explain the fundamental lemma \cite{Willems2005} in a state-space setting. Our goal is to provide a simple and self-contained proof of the result within this context. Consider the linear time-invariant (LTI) system
\begin{subequations}\label{system}
\begin{align} 
\mathbf{x}(t+1) &= A\mathbf{x}(t) + B\mathbf{u}(t) \label{systema} \\
\mathbf{y}(t) &= C\mathbf{x}(t) + D \mathbf{u}(t), \label{systemb}
\end{align}
\end{subequations}
where $\mathbf{x} \in \mathbb{R}^n$ denotes the state, $\mathbf{u} \in \mathbb{R}^m$ is the input and $\mathbf{y} \in \mathbb{R}^p$ is the output. Let $(u_{[0,T-1]},y_{[0,T-1]})$ be a given input/output trajectory\footnote{Throughout this paper, we denote variables such as $\mathbf{u}$ and $\mathbf{y}$ by bold font characters, and specific instances of such variables in normal font, e.g., $u(0),u(1),...$ and $y(0),y(1),...$.} of \eqref{system}. We consider the Hankel matrices of these inputs and outputs, given by:
\begin{equation}
\label{inputoutputHankel}
    \begin{bmatrix}
                \mathcal{H}_{L}(u_{[0,T-1]}) \\
                \mathcal{H}_{L}(y_{[0,T-1]})
    \end{bmatrix} = \begingroup 
    \setlength\arraycolsep{2pt}
    \begin{bmatrix} u(0) & u(1) & \cdots & u(T-L) \\
    \vdots & \vdots & & \vdots \\
    u(L-1) & u(L) & \cdots & u(T-1) \\
    y(0) & y(1) & \cdots & y(T-L) \\
    \vdots & \vdots & & \vdots \\
    y(L-1) & y(L) & \cdots & y(T-1)  
    \end{bmatrix},
    \endgroup
\end{equation}
where $L \geq 1$. Clearly, each column of \eqref{inputoutputHankel} contains a length $L$ input/output trajectory of \eqref{system}. By linearity of the system, every linear combination of the columns of \eqref{inputoutputHankel} is also a trajectory of \eqref{system}. In other words, 
\begin{equation}
\label{parameterization}
    \begin{bmatrix}
        \bar{u}_{[0,L-1]} \\ \bar{y}_{[0,L-1]}
    \end{bmatrix} := \begin{bmatrix}
                \mathcal{H}_{L}(u_{[0,T-1]}) \\
                \mathcal{H}_{L}(y_{[0,T-1]})
    \end{bmatrix} g
\end{equation}
is an input/output trajectory of \eqref{system} for any real vector $g$. 

The powerful crux of Willems \emph{et al.}'s fundamental lemma is that \emph{every} length $L$ input/output trajectory of \eqref{system} can be expressed in terms of $(u_{[0,T-1]},y_{[0,T-1]})$ as in \eqref{parameterization}, assuming that $u_{[0,T-1]}$ is persistently exciting. The result has appeared first in a behavioral context in \cite[Thm. 1]{Willems2005}. In Theorem \ref{theoremWillems}, we will formulate the fundamental lemma for systems of the form \eqref{system}. The theorem consists of two statements. First, under controllability and excitation assumptions, a rank condition on the state and input Hankel matrices \eqref{fullrank} is satisfied.  
Second, under the same conditions, all length $L$ input/output trajectories of \eqref{system} can be written as a linear combination of the columns of the matrix \eqref{inputoutputHankel}. 
\begin{theorem}
\label{theoremWillems}
Consider the system \eqref{system} and assume that the pair $(A,B)$ is controllable. Let $(u_{[0,T-1]},x_{[0,T-1]},y_{[0,T-1]})$ be an input/state/output trajectory of \eqref{system}. Assume that the input $u_{[0,T-1]}$ is persistently exciting of order $n+L$. Then the following statements hold: 
\begin{enumerate}[label=(\roman*),wide, labelindent=0pt]
    \item \label{statement1} The matrix
    \begin{equation}
\label{fullrank}
    \begin{bmatrix} \mathcal{H}_{1}(x_{[0,T-L]}) \\ \mathcal{H}_{L}(u_{[0,T-1]}) \end{bmatrix} = 
    \begingroup 
    \setlength\arraycolsep{2pt}
    \begin{bmatrix} x(0) & x(1) & \cdots & x(T-L) \\ u(0) & u(1) & \cdots & u(T-L) \\ \vdots & \vdots &  & \vdots \\ u(L-1) & u(L) & \cdots & u(T-1) \end{bmatrix}
    \endgroup 
\end{equation}
has full row rank.
    \item \label{statement2} Every length $L$ input/output trajectory of \eqref{system} can be expressed in terms of $u_{[0,T-1]}$ and $y_{[0,T-1]}$ as follows: $(\bar{u}_{[0,L-1]}, \bar{y}_{[0,L-1]})$ is an input/output trajectory of \eqref{system} if and only if 
\begin{equation}
\label{trajectories}
    \begin{bmatrix}
        \bar{u}_{[0,L-1]} \\ \bar{y}_{[0,L-1]}
    \end{bmatrix} = \begin{bmatrix}
                \mathcal{H}_{L}(u_{[0,T-1]}) \\
                \mathcal{H}_{L}(y_{[0,T-1]})
    \end{bmatrix} g,
\end{equation}
for some real vector $g$.
\end{enumerate}
\end{theorem}

Statement \ref{statement1} has appeared first in the original paper by Willems and coworkers, c.f. \cite[Cor. 2(iii)]{Willems2005}. The result is intriguing since a rank condition on \emph{both} input and state matrices can be inposed by injecting a sufficiently exciting input sequence. This rank condition is important from a design perspective and plays a fundamental role in MOESP type subspace algorithms, c.f. \cite[Sec. 3.3]{Verhaegen1992}. Also, in the case that $L = 1$, full row rank of \eqref{fullrank} has been shown to be instrumental for the construction of state feedback controllers from data \cite{DePersis2019}. In our work, statement \ref{statement1} is used to prove the second statement of Theorem \ref{theoremWillems}.
Statement \ref{statement2} is a reformulation of \cite[Thm. 1]{Willems2005}. In what follows, we provide a self-contained and elementary proof of the fundamental lemma in a state-space context.

\begin{IEEEproof}
Statement \ref{statement2} has been proven assuming statement \ref{statement1} in \cite[Lemma 2]{DePersis2019}. It therefore remains to be shown that \eqref{fullrank} has full row rank. Let $\begin{bmatrix} \xi & \eta \end{bmatrix}$ be a vector in the left kernel of \eqref{fullrank}, where $\xi^\top \in \mathbb{R}^n$ and $\eta^\top \in \mathbb{R}^{mL}$. We will first show that $\xi$ and $\eta$ can be used to construct $n+1$ vectors in the left kernel of the ``deeper" Hankel matrix
\begin{equation}
    \label{Hankeln+L}
    \begin{bmatrix} \mathcal{H}_{1}(x_{[0,T-n-L]}) \\ \mathcal{H}_{n+L}(u_{[0,T-1]}) \end{bmatrix}.
\end{equation}
First, by definition of $\xi$ and $\eta$, it is clear that 
$$\begin{bmatrix} \xi & \eta & 0_{nm} \end{bmatrix} \begin{bmatrix} \mathcal{H}_{1}(x_{[0,T-n-L]}) \\ \mathcal{H}_{n+L}(u_{[0,T-1]}) \end{bmatrix} = 0.
$$ 
Next, by the laws of system \eqref{systema} we have
$$
\mathcal{H}_1(x_{[1,T-n-L+1]}) = \begin{bmatrix}
    A & B
\end{bmatrix} \begin{bmatrix}
    \mathcal{H}_1(x_{[0,T-n-L]}) \\
    \mathcal{H}_1(u_{[0,T-n-L]}) 
\end{bmatrix}.
$$
Using this fact, we see that 
\begin{align*}
\begin{bmatrix} \xi A & \xi B & \eta & 0_{(n-1)m} \end{bmatrix}\begin{bmatrix} \mathcal{H}_{1}(x_{[0,T-n-L]}) \\ \mathcal{H}_{n+L}(u_{[0,T-1]}) \end{bmatrix} &= \\
\begin{bmatrix}
    \xi & \eta 
\end{bmatrix}
\begin{bmatrix}
    \mathcal{H}_1(x_{[1,T-n-L+1]}) \\
    \mathcal{H}_L(u_{[1,T-n]})
\end{bmatrix}
&=0,
\end{align*}
where the latter equality holds by definition of $\xi$ and $\eta$. Now, by repeatedly exploiting the laws of \eqref{systema} and using the same arguments we find that the $n+1$ vectors
\begin{equation}
\label{depvectors}
\begin{aligned}
      &w_0 := \begin{bmatrix} \xi & \eta & 0_{nm} \end{bmatrix} \\
      &w_1 := \begin{bmatrix} \xi A & \xi B & \eta & 0_{(n-1)m} \end{bmatrix} \\
      &w_2 := \begin{bmatrix} \xi A^2 & \xi AB & \xi B & \eta & 0_{(n-2)m} \end{bmatrix} \\
      &\: \vdots  \\
      &w_n := \begin{bmatrix} \xi A^n & \xi A^{n-1} B & \cdots & \xi B & \eta \end{bmatrix}
\end{aligned}
\end{equation}
are all contained in the left kernel of the matrix \eqref{Hankeln+L}. By persistency of excitation, $\mathcal{H}_{n+L}(u_{[0,T-1]})$ has full row rank, and hence the left kernel of \eqref{Hankeln+L} has dimension at most $n$. Therefore, the $n+1$ vectors in \eqref{depvectors} are linearly dependent. We claim that this implies $\eta = 0$. To prove this claim, partition $\eta = \begin{bmatrix}
    \eta_1 & \eta_2 & \cdots & \eta_L
\end{bmatrix}$, where $\eta_1^\top,\eta_2^\top,\dots,\eta_L^\top \in \mathbb{R}^m$. Since the last $m$ entries of the vectors $w_0,w_1,\dots,w_{n-1}$ are zero, the linear dependence of the vectors \eqref{depvectors} implies $\eta_L = 0$ by inspection of $w_n$. We substitute this equation in $\eta$ and conclude that the last $2m$ entries of $w_0,w_1,\dots,w_{n-1}$ are zero. As such, also $\eta_{L-1} = 0$. We can proceed with these substitutions to show that $\eta_1 = \eta_2 = \cdots \eta_L = 0$, i.e., $\eta = 0$. Next, by Cayley-Hamilton theorem, $\sum_{i = 0}^{n} \alpha_i A^i = 0$ where $\alpha_i \in \mathbb{R}$ for all $i = 0,1,\dots,n$, and $\alpha_n = 1$. Define the linear combination $v := \sum_{i=0}^n \alpha_i w_i$. By \eqref{depvectors} and by substitution of $\eta = 0$, the vector $v$ is equal to
$$\begingroup\setlength\arraycolsep{3pt}\begin{bmatrix} 0_n & \sum_{i=1}^n \alpha_i \xi A^{i-1} B & \sum_{i=2}^n \alpha_i \xi A^{i-2} B & \cdots & \alpha_n \xi B & 0_{mL} \end{bmatrix}.\endgroup
$$
This implies that the vector 
$$
\begin{bmatrix} \sum_{i=1}^n \alpha_i \xi A^{i-1} B & \sum_{i=2}^n \alpha_i \xi A^{i-2} B & \cdots & \alpha_n \xi B \end{bmatrix}
$$
is contained in the left kernel of $\mathcal{H}_{n}(u_{[0,T-L-1]})$, which is zero by persistency of excitation. In other words, 
\begin{align*}
0 &= \alpha_1 \xi B + \cdots + \alpha_n \xi A^{n-1} B \\
0 &= \alpha_2 \xi B + \cdots + \alpha_n \xi A^{n-2} B \\
&\:\:\vdots \\
0 &= \alpha_{n-1} \xi B + \alpha_n \xi AB \\
0 &= \alpha_n \xi B.
\end{align*}
Since $\alpha_n = 1$ it follows from the last equation that $\xi B = 0$. Substitution in the second to last equation then results in $\xi AB = 0$. We continue by backward substitution to obtain $\xi B = \xi AB = \cdots = \xi A^{n-1}B = 0$. Controllability of $(A,B)$ hence results in $\xi = 0$. We therefore conclude that \eqref{fullrank} has full row rank, which proves the theorem. 
\end{IEEEproof}

\section{Extension of Willems \emph{et al.}'s lemma to multiple trajectories}
\label{sectionmultipledata}

In this section we propose an extension of the fundamental lemma that is applicable to the case in which \emph{multiple} system trajectories are given. Our approach will require the notion of \emph{collective} persistency of excitation. 

\begin{definition}
Consider the input sequences $u^i_{[0,T_i-1]}$ for $i = 1,2,\dots,q$, where $q$ is the number of data sets. Let $k$ be a positive integer such that $k \leq T_i$ for all $i$. The input sequences $u^i_{[0,T_i-1]}$ for $i = 1,2,\dots,q$ are called \emph{collectively persistently exciting} of order $k$ if the mosaic-Hankel matrix
\begin{equation}
\label{mosaicHankel}
    \begingroup
    \setlength\arraycolsep{2pt}
        \begin{bmatrix}
            \mathcal{H}_{k}(u^1_{[0,T_1-1]}) & \mathcal{H}_{k}(u^2_{[0,T_2-1]}) & \cdots & \mathcal{H}_{k}(u^q_{[0,T_q-1]})
        \end{bmatrix}
        \endgroup
\end{equation}
has full row rank. 
\end{definition}

Collective persistency of excitation is more flexible than the persistency of excitation of a single input sequence. Indeed, for the input sequences $u^i_{[0,T_i]}$ to be collectively persistently exciting, it is sufficient that at least one of them is persistently exciting. However, this is clearly not necessary: the sequences $u^i_{[0,T_i]}$ may be collectively persistently exciting even when 
none of the individual input sequences is persistently exciting. The added flexibility of collective persistency of excitation is also apparent from the \emph{length} of the input sequences. Indeed, a single $u_{[0,T-1]}$ can only be persistently exciting of order $k$ if $T \geq k(m+1)-1$. In comparison, for collective persistency of excitation of order $k$ it is necessary that $\sum_{i = 1}^q T_i \geq k(m+q)-q$. This means that collective persistency of excitation can be achieved by input sequences having length $T_i$ as short as $k$, assuming the number of data sets $q$ is sufficiently large. In the next theorem we extend the fundamental lemma to the case of multiple data sets.

\begin{theorem}
\label{theoremmultipledata}
Consider system \eqref{system} and assume that the pair $(A,B)$ is controllable. Let 
$(u^i_{[0,T_i-1]},x^i_{[0,T_i-1]},y^i_{[0,T_i-1]})$
be an input/state/output trajectory of \eqref{system} for 
$i = 1,2,\dots,q$. Assume that the inputs $u^i_{[0,T_i-1]}$ are collectively persistently exciting of order $n+L$. Then the following statements hold: 
\begin{enumerate}[label=(\roman*),wide, labelindent=0pt]
    \item \label{statement1md} 
    The matrix
\begin{equation}
\label{jointfullrank}
\begingroup
    \setlength\arraycolsep{2pt}
    \begin{bmatrix} \mathcal{H}_{1}(x^1_{[0,T_1-L]}) & \mathcal{H}_{1}(x^2_{[0,T_2-L]}) & \cdots & \mathcal{H}_{1}(x^q_{[0,T_q-L]}) \\ \mathcal{H}_{L}(u^1_{[0,T_1-1]}) & \mathcal{H}_{L}(u^2_{[0,T_2-1]}) & \cdots & \mathcal{H}_{L}(u^q_{[0,T_q-1]}) \end{bmatrix} \endgroup
\end{equation}
has full row rank. 
 \item \label{statement2md} 
Every length $L$ input/output trajectory of \eqref{system} can be expressed in terms of $u^i_{[0,T_i-1]}$ and $y^i_{[0,T_i-1]}$ ($i = 1,2,\dots,q$) as follows: $(\bar{u}_{[0,L-1]}$, $\bar{y}_{[0,L-1]})$ is an input/output trajectory of \eqref{system} if and only if
\begin{align}
\label{jointtrajectories}
\begin{bmatrix} \bar{u}_{[0,L-1]} \\ \bar{y}_{[0,L-1]}
    \end{bmatrix} = 
\begingroup
    \setlength\arraycolsep{2pt}\begin{bmatrix}
                \mathcal{H}_{L}(u^1_{[0,T_1-1]}) & 
                \cdots & \mathcal{H}_{L}(u^q_{[0,T_q-1]}) \\
                \mathcal{H}_{L}(y^1_{[0,T_1-1]}) & 
                \cdots & \mathcal{H}_{L}(y^q_{[0,T_q-1]})
    \end{bmatrix}\endgroup g,
\end{align}
for some real vector $g$.
\end{enumerate}
\end{theorem}

Note that if $q = 1$ and $T_1 = T$ we deal with a single experiment, and in this case Theorem \ref{theoremmultipledata} recovers Theorem \ref{theoremWillems}.

 \begin{IEEEproof}
We first prove that \eqref{jointfullrank} has full row rank. Let $\begin{bmatrix} \xi & \eta \end{bmatrix}$ be a vector in the left kernel of \eqref{jointfullrank}, where $\xi^\top \in \mathbb{R}^n$ and $\eta^\top \in \mathbb{R}^{mL}$. By exploiting the laws of the system \eqref{systema} we see that the vectors
\begin{equation}
\label{jointdepvectors}
\begin{aligned}
      &w_0 := \begin{bmatrix} \xi & \eta & 0_{nm} \end{bmatrix} \\
      &w_1 := \begin{bmatrix} \xi A & \xi B & \eta & 0_{(n-1)m} \end{bmatrix} \\
      &w_2 := \begin{bmatrix} \xi A^2 & \xi AB & \xi B & \eta & 0_{(n-2)m} \end{bmatrix} \\
      &\: \vdots  \\
      &w_n := \begin{bmatrix} \xi A^n & \xi A^{n-1} B & \cdots & \xi B & \eta \end{bmatrix}
\end{aligned}
\end{equation}
are contained in the left kernel of the matrix 
\begin{equation}
    \label{jointHankeln+L}
    \begin{bmatrix} \mathcal{H}_{1}(x^1_{[0,T_1-n-L]}) & \cdots & \mathcal{H}_{1}(x^q_{[0,T_q-n-L]}) \\ \mathcal{H}_{n+L}(u^1_{[0,T_1-1]}) & \cdots & \mathcal{H}_{n+L}(u^q_{[0,T_q-1]}) \end{bmatrix}.
\end{equation}
By the persistency of excitation assumption, the matrix
$$
\begin{bmatrix} \mathcal{H}_{n+L}(u^1_{[0,T_1-1]}) & \cdots & \mathcal{H}_{n+L}(u^q_{[0,T_q-1]}) \end{bmatrix}
$$
has full row rank, and hence the left kernel of \eqref{jointHankeln+L} has dimension at most $n$. Therefore, the $n+1$ vectors in \eqref{jointdepvectors} are linearly dependent. This yields $\eta = 0$ following the same argument as in the proof of Theorem \ref{theoremWillems}. Next, by Cayley-Hamilton theorem, $\sum_{i = 0}^{n} \alpha_i A^i = 0$ where $\alpha_i \in \mathbb{R}$ for $i = 0,1,\dots,n$ and $\alpha_n = 1$. We define the linear combination $v := \sum_{i=0}^n \alpha_i w_i$. Clearly, the vector $v$ is equal to
$$\begingroup\setlength\arraycolsep{3pt}\begin{bmatrix} 0_n & \sum_{i=1}^n \alpha_i \xi A^{i-1} B & \sum_{i=2}^n \alpha_i \xi A^{i-2} B & \cdots & \alpha_n \xi B & 0_{mL} \end{bmatrix}.\endgroup
$$
Hence, the vector 
$$
\begin{bmatrix} \sum_{i=1}^n \alpha_i \xi A^{i-1} B & \sum_{i=2}^n \alpha_i \xi A^{i-2} B & \cdots & \alpha_n \xi B \end{bmatrix}
$$
is contained in the left kernel of
$$
\begin{bmatrix} \mathcal{H}_{n}(u^1_{[0,T_1-L-1]}) & \cdots & \mathcal{H}_{n}(u^q_{[0,T_q-L-1]}) \end{bmatrix},
$$
which is zero by collective persistency of excitation. Following the same steps as in the proof of Theorem \ref{theoremWillems} we conclude by backward substitution that $\xi B = \xi AB = \cdots = \xi A^{n-1}B = 0$. By controllability of $(A,B)$ we have $\xi = 0$, proving statement \ref{statement1md}. 

Next, we prove statement \ref{statement2}. 
Let $\bar{u}_{[0,L-1]}$ and $\bar{y}_{[0,L-1]}$ be vectors such that \eqref{jointtrajectories} is satisfied for some $g$. Then $$\begin{bmatrix} \bar{u}_{[0,L-1]} \\ \bar{y}_{[0,L-1]} \end{bmatrix}$$ is a linear combination of length $L$ trajectories of \eqref{system} and hence, by linearity, itself an input/output trajectory of \eqref{system}. Conversely, let $(\bar{u}_{[0,L-1]}$, $\bar{y}_{[0,L-1]})$ be an input/output trajectory of \eqref{system} and denote by $\bar{x}_0$ a corresponding initial state at time $0$. We have the relation
\begin{equation}
\label{OL}
    \begin{bmatrix}
        \bar{u}_{[0,L-1]} \\ \bar{y}_{[0,L-1]}
    \end{bmatrix} = \begin{bmatrix}
                0 & I \\ \mathcal{O}_L & \mathcal{T}_L
    \end{bmatrix} \begin{bmatrix}
        \bar{x}_0 \\ \bar{u}_{[0,L-1]} 
    \end{bmatrix},
\end{equation}
where $\mathcal{T}_L$ and $\mathcal{O}_L$ are defined as
\begin{align}
    \label{toeplitz}
\mathcal{T}_L &:= \begingroup\setlength\arraycolsep{2pt}\begin{bmatrix}
            D & 0 & 0 & \cdots & 0 \\
            CB & D & 0 & \cdots & 0 \\
            CAB & CB & D & \cdots & 0 \\
            \vdots & \vdots & \vdots & \ddots & \vdots \\
            CA^{L-2}B & CA^{L-3}B & CA^{L-4}B & \cdots & D
\end{bmatrix},
\endgroup 
\\
\label{markov}
\mathcal{O}_L &:= \begingroup\setlength\arraycolsep{2pt}\begin{bmatrix}
            C^\top & (CA)^\top & (CA^2)^\top & \cdots & (CA^{L-1})^\top
\end{bmatrix}^\top.
\endgroup
\end{align}
Since \eqref{jointfullrank} has full row rank, there exists a vector $g$ such that 
$$
\begin{bmatrix}
        \bar{x}_0 \\ \bar{u}_{[0,L-1]} 
    \end{bmatrix} =  \begin{bmatrix} \mathcal{H}_{1}(x^1_{[0,T_1-L]}) & \cdots & \mathcal{H}_{1}(x^q_{[0,T_q-L]}) \\ \mathcal{H}_{L}(u^1_{[0,T_1-1]}) & \cdots & \mathcal{H}_{L}(u^q_{[0,T_q-1]}) \end{bmatrix} g.
$$
Substitution of the latter expression into \eqref{OL} and using the fact that
$$
\begin{bmatrix}
                0 & I \\ \mathcal{O}_L & \mathcal{T}_L
    \end{bmatrix} \begin{bmatrix} \mathcal{H}_{1}(x^i_{[0,T_i-L]}) \\ \mathcal{H}_{L}(u^i_{[0,T_i-1]}) \end{bmatrix} = \begin{bmatrix}
                \mathcal{H}_{L}(u^i_{[0,T_i-1]}) \\
                \mathcal{H}_{L}(y^i_{[0,T_i-1]})
    \end{bmatrix}
$$
for all $i = 1,2,\dots,q$ yields \eqref{jointtrajectories}, as desired.
\end{IEEEproof}

\section{Examples of application}
\label{sectionexamples}

\subsection{Identification with missing data samples}
In this section we treat an example in which we want to identify a system model from a measured trajectory with missing data samples. System identification from trajectories with missing data has been studied in the papers \cite{Markovsky2013,Markovsky2016,Markovsky2017}. As we will see, it is also possible to apply Theorem \ref{theoremmultipledata}\ref{statement2md} in this context.

Suppose that we have access to the following, partially corrupted, 
input/output trajectory of length $T = 20$: 
\begin{align*}
\begin{array}{ |c|r|r|r|r|r|r|r|r|r|r| } 
 t    & 0 & 1 & 2 & 3 & 4 & 5 & 6 & 7 & 8 & 9 \\ 
 u(t) & 1 & 0 & 2 & -1 & 0 & \boldsymbol \times & 1 & 1 & -1 & -5 \\ 
 y(t) & 3 & 3 & 7 & 6 & 11 & \boldsymbol \times & 18 & 21 & 23 & 24 \\ 
\end{array}
\end{align*}
\begin{align*}
\begin{array}{ |c|r|r|r|r|r|r|r|r|r|r| } 
 t    & 10 & 11 & 12 & 13 & 14 & 15 & 16 & 17 & 18 & 19 \\ 
 u(t) & 0 & -1 & \boldsymbol \times & 1 & -6 & 2 & -2 & 0 & 1 & \boldsymbol \times \\ 
 y(t) & 33 & 31 & \boldsymbol \times & 30 & 20 & 26 & 14 & 10 & 3 & \boldsymbol \times \\ 
\end{array}
\end{align*}

The data are generated by a minimal LTI system of (unknown) state-space dimension $n=2$. Note that some of the samples are \emph{missing}, which we indicate by $\boldsymbol \times$. Our goal is to identify an LTI system that is compatible with the observed data. 

In this problem, we have
access to three input/output system trajectories, namely $(u_{[0,4]},y_{[0,4]})$, $(u_{[6,11]},y_{[6,11]})$ and $(u_{[13,18]},y_{[13,18]})$. It is not difficult to verify that the input sequences $u_{[0,4]}$, $u_{[6,11]}$ and $u_{[13,18]}$ are collectively persistently exciting of order $5$. It can be easily verified that no LTI system of dimension $1$ can explain the data. Thus we consider LTI systems of dimension $2$. Since the inputs are collectively persistently exciting of order $5$, and since the data-generating system has dimension $n=2$, by Theorem \ref{theoremmultipledata}\ref{statement2md} every length $L=3$
input/output trajectory of the system can be written as linear combination of the columns of 
\begin{equation}
    \label{datahankel}
    \mathcal{D} :=
\begingroup
    \setlength\arraycolsep{2pt}\begin{bmatrix}
                \mathcal{H}_3(u_{[0,4]}) & \mathcal{H}_3(u_{[6,11]}) & \mathcal{H}_3(u_{[13,18]}) \\
                \mathcal{H}_3(y_{[0,4]}) & \mathcal{H}_3(y_{[6,11]}) & \mathcal{H}_3(y_{[13,18]})
    \end{bmatrix}. \endgroup
\end{equation}
We exploit this result by computing, as a function of $\mathcal{D}$, the length $7$ system trajectory 
\begin{align}
\bar{u}_{[-2,4]} &= \begin{bmatrix}
                0 & 0 & 1 & 0 & 0 & 0 & 0 
\end{bmatrix}^\top \\
\label{ybar}
\bar{y}_{[-2,4]} &= \begin{bmatrix}
                0 & 0 & \textbf{?} & \textbf{?} & \textbf{?} & \textbf{?} & \textbf{?}
\end{bmatrix}^\top, 
\end{align}
where question marks denote to-be-computed values. The idea is as follows: 
if the ``past" inputs $\bar{u}(-2),\bar{u}(-1)$ and ``past" outputs $\bar{y}(-2),\bar{y}(-1)$ are zero, the state $\bar{x}(0) \in \mathbb{R}^2$ corresponding to $(\bar{u}_{[-2,4]},\bar{y}_{[-2,4]})$ is unique, and equal to zero. 
This means that $\bar{u}_{[0,4]}$ is an impulse, applied to a system of the form \eqref{system} 
with zero initial state. Consequently, the output $\bar{y}_{[0,4]}$ simply consists of the first 
Markov parameters of \eqref{system}, that is, $\bar{y}_{[0,4]} = \begin{bmatrix}
                D & CB & CAB & CA^2B & CA^3B
\end{bmatrix}$. From these Markov parameters it is straightforward to 
compute a state-space realization, e.g., using the Ho-Kalman algorithm
\cite[Section 3.4.4]{Verhaegen2007}.
 
Therefore, our remaining task is to compute $\bar{y}_{[0,4]}$. Inspired by \cite{Markovsky2005b}, we will compute this trajectory iteratively by computing multiple length 3 trajectories as linear combinations of the columns of \eqref{datahankel}. To begin with, 
we compute the first unknown in \eqref{ybar}, which is $\bar{y}(0)$. To do so, we have to solve the system of linear equations\footnote{
Note that the the solution $g$ is not unique in general, but $\bar{y}(0)$ \emph{is} unique. The reason is that the initial state $\bar{x}(0) = 0$ is uniquely specified by the ``past" inputs $\bar{u}(-2),\bar{u}(-1)$ and outputs $\bar{y}(-2),\bar{y}(-1)$. In turn, the initial state $\bar{x}(0)$ and input $\bar{u}(0)$ uniquely specify the output $\bar{y}(0)$. Also see \cite[Prop. 1]{Markovsky2008}.
} 
\begin{equation}
\label{Dg}
    \mathcal{D} g = \begin{bmatrix}
                    0 & 0 & 1 & 0 & 0 & \bar{y}(0)
    \end{bmatrix}^\top 
\end{equation}
in the unknowns $g$ and $\bar{y}(0)$. One possible approach \cite[Alg. 1]{Markovsky2008} is to obtain a solution $\bar{g}$ to the first five linear equations in \eqref{Dg}. Subsequently, $\bar{y}(0)$ is obtained by multiplication of the last row of $\mathcal{D}$ with $\bar{g}$. We do this to find $\bar{y}(0) = 1$. Next, to find $\bar{y}(1)$ we complete the length 3 trajectory $(\bar{u}_{[-1,1]},\bar{y}_{[-1,1]})$ by solving the system of equations
\begin{equation*}
    \mathcal{D} g = \begin{bmatrix}
                    0 & 1 & 0 & 0 & 1 & \bar{y}(1)
    \end{bmatrix}^\top,
\end{equation*}
which results in $\bar{y}(1) = 0$. Repeating this process, we obtain $\bar{y}(2) = 1$, $\bar{y}(3) = 2$ and $\bar{y}(4) = 3$, meaning that
$$
D = 1, \:\: CB = 0, \:\: CAB = 1, \:\: CA^2 B = 2, \:\: CA^3 B = 3.
$$
Finally, it is not difficult to obtain a state-space realization of these Markov parameters as
\begin{align*} 
A = \begin{bmatrix}
                1 & 0 \\ 1 & 1
\end{bmatrix}, \: B = \begin{bmatrix}
                1 \\ 0
\end{bmatrix}, \: C = \begin{bmatrix}
                0 & 1
\end{bmatrix}, \: D = 1.
\end{align*}

The approach outlined in this section is generally also applicable in the case that multiple consecutive data samples are missing. Even in the case that the number of consecutive missing samples is \emph{unknown}, we can apply Theorem \ref{theoremmultipledata} to the partial trajectories. Note that we require a sufficient number of partial trajectories of length at least 5 to guarantee collective persistency of excitation of order $5$. In the case of missing data with larger frequency, it may still be possible to identify the system by computation of the left kernels of submatrices of the Hankel matrix \cite{Markovsky2016}.

\subsection{Data-driven LQR of an unstable system}

Consider the unstable batch reactor system \cite{Walsh2001}, 
which we have discretized using a sampling time of $0.5$s to obtain a system of the form \eqref{systema} with
$$ 
\resizebox{0.99\hsize}{!}{%
$A = \begin{bmatrix}
                2.622 & 0.320 & 1.834 & -1.066 \\
                -0.238 & 0.187 & -0.136 & 0.202 \\
                0.161 & 0.789 & 0.286 & 0.606 \\
                -0.104 & 0.764 & 0.089 & 0.736 
\end{bmatrix}, \:\: B = \begin{bmatrix}
                0.465 & -1.550 \\
                1.314 & 0.085 \\
                2.055 & -0.673 \\
                2.023 & -0.160 
\end{bmatrix}$.
}
$$
The goal of this example is the data-based design of an optimal control input $u^*$ that minimizes the cost functional 
$$
J := \sum_{t = 0}^{\infty} x^\top(t)Qx(t) + u^\top(t) R u(t) .
$$
under the zero endpoint constraint $\lim_{t \to \infty} x(t) = 0$. Here $Q$ and $R$ are state and input weight matrices, respectively. Under standard assumptions on $A$, $B$, $Q$ and $R$ \cite[Thm. 23]{vanWaarde2020}, 
the optimal input exists, is unique, and 
{is generated by the feedback law $u^* = Kx$, where
$$
K = -(R+B^\top P^+ B)^{-1} B^\top P^+ A
$$
and where $P^+$ is the largest real symmetric solution to the algebraic Riccati equation 
$$
P = A^\top P A - A^\top P B(R+B^\top PB)^{-1} B^\top PA +Q.
$$

In \cite[Thm. 4]{DePersis2019} an attractive design procedure is 
introduced to obtain $K$ directly from input/state data. 
 The idea is to inject an input sequence $u_{[0,T-1]}$ that is persistently exciting of order $n+1$ 
 such that the matrix\footnote{Note that
 $X_- := \mathcal{H}_1(x_{[0,T-1]})$ and $U_- := \mathcal{H}_1(u_{[0,T-1]})$.
 }
\begin{equation} 
\label{xuexample}
\begin{bmatrix}
X_- \\ \hline U_-
\end{bmatrix} := 
\begin{bmatrix}
    x(0) & x(1) & \cdots & x(T-1) \\ \hline
    u(0) & u(1) & \cdots & u(T-1)
\end{bmatrix}
\end{equation}
has full row rank by Theorem \ref{theoremWillems}\ref{statement1}. 
Subsequently, $K$ is found by solving a semidefinite program involving 
the data $x_{[0,T]}$ and $u_{[0,T-1]}$ alone; see \cite[Eq. 27]{DePersis2019}. Later on, it was shown \cite[Thm. 26]{vanWaarde2020} that full row rank of \eqref{xuexample} is actually also \emph{necessary} for obtaining $K$ from input/state data. In addition, another semidefinite program was introduced \cite[Thm. 29]{vanWaarde2020} to obtain $P^+$ and $K$ from input/state data.  
Both semidefinite programs of \cite{DePersis2019} and \cite{vanWaarde2020} are applicable to this example, but we will follow the method of \cite{vanWaarde2020} since it involves less decision variables, c.f. \cite[Remark 31]{vanWaarde2020}. We will compare the approach based on 
a \emph{single} measured trajectory of the system with the one
based on \emph{multiple} trajectories. In both the approaches, 
we take $Q$ and $R$ as the identity matrices of appropriate dimensions.

First, we compute $K$ on the basis of a \emph{single} measured trajectory of \eqref{systema}. 
We choose a random initial state and random input sequence of length $T = 20$, 
generated using the Matlab command \texttt{rand}. This input is persistently exciting 
of order $5$. 
Finally, we let $X_-$ and $U_-$ as in \eqref{xuexample}, and define
$X_+ := \mathcal{H}_1(x_{[1,T]})$. By \cite[Thm. 29]{vanWaarde2020}, 
the largest solution $P^+$ to the algebraic Riccati equation is the unique solution to the optimization problem

\begin{equation}
\label{opt}
\begin{aligned}
    \text{maximize } &\tr P \\
    \text{subject to } &P = P^\top \geq 0 \text{ and } \mathcal{L}(P) \leq 0,
\end{aligned}
\end{equation}
where $\mathcal{L}(P) := X_-^\top P X_- -X_+^\top P X_+ -X_-^\top Q X_- - U_-^\top R U_-$.
We use Yalmip with Sedumi 1.3 as LMI solver. 
Because of the large magnitude of the data samples (reaching $||x(19)|| = 1.049 \cdot 10^{8}$), the solver runs into 
numerical problems and returns a matrix $P_{\sing}$ that does not resemble $P^+$. 
In fact, comparing $P_{\sing}$ with the ``true" matrix $P^+$ obtained via the (model-based) 
Matlab command \texttt{dare}, we see

\begin{align*} 
P_{\sing} &= \begin{bmatrix}
    0.002 & 0.013 & -0.005 & 0.015 \\
    0.013 & 0.075 & 0.017 & 0.067 \\
    -0.005 & 0.017 & 0.823 & 0.066 \\
    0.015 & 0.067 & 0.066 & 0.010 
\end{bmatrix} \\
P^+ &= \begin{bmatrix}
    3.604 & 0.049 & 1.762 & -1.306 \\
    0.049 & 1.170 & 0.072 & 0.142 \\
    1.762 & 0.072 & 2.202 & -0.845 \\
    -1.306 & 0.142 & -0.845 & 1.823
\end{bmatrix}.
\end{align*}
To overcome this problem, we next consider \emph{multiple} short 
experiments, demonstrating the effectiveness of this second approach. 
We collect $q = 5$ data sets of length $T_i = 6$ for $i = 1,2,3,4,5$. 
The input sequences $u^i_{[0,T_i-1]}$ of these sets are again chosen randomly, 
and are verified to be collectively persistently exciting of order $5$. Similar as before, we use 
the notation $X^i_- := \mathcal{H}_1(x^i_{[0,T_i-1]})$, 
$X^i_+ := \mathcal{H}_1(x^i_{[1,T_i]})$ and $U^i_- := \mathcal{H}_1(u^i_{[0,T_i-1]})$ for all $i$. 
In addition, we concatenate these data matrices and define
\begin{align*}
    X_- &:= \begin{bmatrix}
        X_-^1 & X_-^2 & \cdots & X_-^5 
    \end{bmatrix} \\
    X_+ &:= \begin{bmatrix}
        X_+^1 & X_+^2 & \cdots & X_+^5 
    \end{bmatrix} \\
    U_- &:= \begin{bmatrix}
        U_-^1 & U_-^2 & \cdots & U_-^5 
    \end{bmatrix}.
\end{align*}
With these data matrices, we solve again \eqref{opt}.
This result in the solution $P_{\mult}$  
with $||P_{\mult} - P^+|| = 7.849 \cdot 10^{-10}$. 
Next, we continue the design procedure of \cite[Thm. 29]{vanWaarde2020} by computing a right inverse $X_-^\dagger$ of $X_-$ such that $\mathcal{L}(P_{\mult})X_-^\dagger = 0$. 
The optimal control gain is then computed as 
$$
K_{\mult} := U_-X_-^\dagger = \begin{bmatrix}
    0.163 & -0.292 & 0.046 & -0.328 \\
    1.418 & 0.116 & 0.984 & -0.625
\end{bmatrix}.
$$
The error between between $K_{\mult}$ and the true optimal gain $K$ obtained via the command \texttt{dare} is small. In fact, we have $||K_{\mult} - K|| = 7.083 \cdot 10^{-11}$. The closed-loop matrix $A + BK_{\mult}$ is stable and its spectral radius is $0.188$. 

The approach that uses multiple trajectories overall requires more samples than the one using a single trajectory. Indeed, as explained in Section \ref{sectionmultipledata}, a necessary condition for collective persistency of excitation of order $k$ is that 
$$
\sum_{i = 1}^q T_i \geq k(m+q) -q.
$$
This means that $\sum_{i = 1}^5 T_i \geq 30$ in our example. In comparison, a necessary condition for persistency of excitation of order $5$ of a single trajectory is $T \geq 14$. Nonetheless, as shown in this example, the use of multiple short trajectories enables the accurate computation of feedback gains even for unstable systems while this may be problematic when using a single long trajectory.

\section{Conclusions}
\label{sectionconclusions}
Willems \emph{et al.}'s fundamental lemma is a beautiful result that asserts that all trajectories of a linear system can be parameterized by a single, persistently exciting one. 
In this paper we have extended the fundamental lemma to the scenario where multiple trajectories are given instead of a single one. To this end, we have introduced a notion of collective persistency of excitation. Subsequently, we have shown that all trajectories of a linear system can be parameterized by a finite number of them, assuming these are collectively persistently exciting. We have shown that this result enables the identification of linear systems from data sets with missing data samples. We have also shown that the result can be used to construct controllers of unstable systems from multiple measured trajectories, even when this is not possible from a single trajectory.

\bibliographystyle{IEEEtran}
\bibliography{references}

\end{document}